# THE DIMENSIONS OF HAUSDORFF AND MENDÈS FRANCE.
# A COMPARATIVE STUDY

**R. Hansen and M. Piacquadio.**

Dpto. de Matemática, Fac. de Ingeniería, Universidad de Buenos Aires.

**Abstract**

This paper contains a comparative study of two families of simple curves drawn in the plane. On the one hand, we have the fractal curves on the unit interval, with self-similar structure, which have associated a Hausdorff dimension. On the other hand, we have the opposite: a class of locally rectifiable unbounded curves, which have another "fractional dimension" defined by M. Mendès France. We propose a geometrical constructive process that will allow us to obtain —as the limit of a sequence of polygonal curves— one curve of the first family, by contractive transformations; and another of the second family, by expansive transformations. Thanks to this process of linking curves from both families, we are able to compare their dimensions —our aim in this work—, and to obtain interesting results such as the equality of the latter in the case of strict self-similarity.

*...The reader may feel surprised that there is no mention of Benoît Mandelbrot in these notes. His objects are fractals, i.e., locally irregular. Mine, on the contrary are locally smooth. The curves I discuss are locally rectifiable. My topic could be thought of "anti-Mandelbrotian" within "Mandelbrotmania". I was, I am, and I hope to remain influenced by B. Mandelbrot.* [1]

<div style="text-align:right">Michel Mendès France.</div>

## 1. INTRODUCTION

In this paper we will study two families of non-intersecting planar curves. The first of these families, which we call $\mathcal{F}_H$, is composed of fractal curves with a self-similar structure defined by $N$ contractions of ratios $a_1, a_2, ..., a_N$ ($0 < a_i < 1$, $i = 1, ..., N$) and satisfying the *Closed-set criterion* ([2]): this criterion guarantees their non-intersecting character. We are interested in the Hausdorff dimension $d_H$ of these curves, which satisfies the equation:

$$\sum_{i=1}^{N} a_i^{d_H} = 1 \ .$$

The second family, that we will call $\mathcal{F}_{MF}$, contains locally rectifiable unbounded curves, i.e. any arc of the curve $\Gamma$ has finite length, and we are interested in the dimension $d_{MF}$ introduced by Mendès France in [1], for this type of curve.

The two families have no curve in common; moreover, they have absolutely different geometrical properties; however, we will see that there exists a geometrical constructive process that will allow us to link curves of both families, and thus, to be able to compare their respective dimensions —such is the aim in this work.

## 2. THE DIMENSION OF MENDÈS FRANCE OF A CURVE $\Gamma$

For a curve $\Gamma$ belonging to $\mathcal{F}_{MF}$, we fix an origin and we consider the first portion $\Gamma_L$ of $\Gamma$ of length $L$. Let $\varepsilon > 0$ be given, and let us consider the set:

$$\Gamma_L(\varepsilon) \ = \ \{P \in I\!\!R^2 / \mathrm{dist}(P, \Gamma_L) < \varepsilon/2\} \ .$$

This set is also known as the $\varepsilon$-*Minkowski sausage* of $\Gamma_L$. Let $C_L$ be the length of the boundary of the convex hull of $\Gamma_L$. Then, the dimension of Mendès France of a curve $\Gamma$ is, by definition:

$$\dim_{MF}(\Gamma) \ = \ \lim_{\varepsilon \searrow 0} \liminf_{L \nearrow \infty} \frac{\log A(\Gamma_L(\varepsilon))}{\log C_L} \ ,$$

where $A(\Gamma_L(\varepsilon))$ denotes the area of $\Gamma(\varepsilon)$. There is a remark in [1] that shows that the last formula **does not** depend on $\varepsilon$, so we will either take it away, or replace it by a suitable value in order to make calculations; hence we can write:

$$\dim_{MF}(\Gamma) \;=\; \liminf_{L \nearrow \infty} \frac{\log A(\Gamma_L(\varepsilon))}{\log C_L} \;\;.$$

This remark is very important, because, intuitively, it says that it doesn't matter how "wide" the $\varepsilon$-Minkowski sausage is, but how the sausage "fills up" the plane according to the development of $\Gamma_L$ when $L$ grows. Therefore, we are dealing with a type of dimension which does not look at the curve with a "zoom lens" —as the Hausdorff dimension does; on the contrary, this dimension looks from afar at the behavior of the curve when its length tends to infinity.

To illustrate this idea, let us consider two well known curves. First, the *Archimedean spiral* of step equal to $r$ (Fig.1). When the length $L$ tends to infinity, we have to step away from the plane again and again to observe its behavior, because its convex hull also grows. And if we continue moving away, soon we won't be able to distinguish the step $r$, and we will see the spiral as filling up $I\!\!R^2$ completely. Let us now take $\varepsilon = r$, and let us consider the corresponding $\varepsilon$-Minkowski sausage of $\Gamma$, $\Gamma(r)$. We can easily see that it covers all the plane; then the dimension of Mendès France of $\Gamma$ is 2. Instead, if we consider a *logarithmic spiral* (Fig. 2), it doesn't matter which is the value of $\varepsilon$ chosen; no matter how far away we are from $I\!\!R^2$, we always see an arc —the same arc— of the curve $\Gamma$, which has dimension of Hausdorff equal to 1, and its $\Gamma(\varepsilon)$ will always appear equally "thin". It comes therefore, as no surprise that the dimension of Mendès France of this curve is unity.

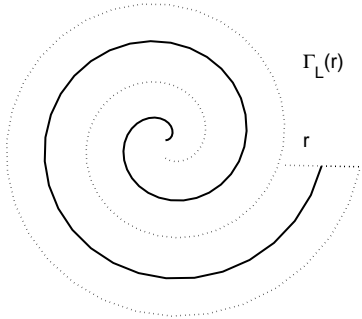
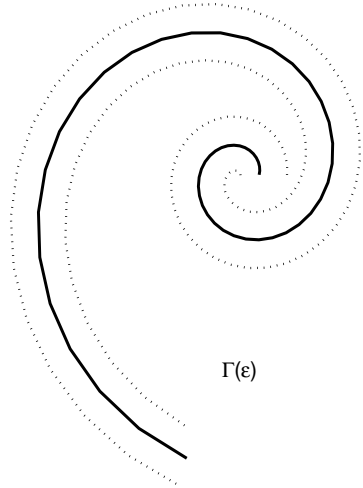

**Figure 1**        **Figure 2**

3. THE STRICT SELF-SIMILAR CURVES

The Hausdorff dimension of a fractal is not, in general, easy to compute, unless the fractal has, for example, some self-similar structure.
Among these cases we have the fractal curves on the unit square interval, with strict self-similarity, like the *Von Koch curve* (Fig. 3).

The process by which we obtain such a curve, consists of replacing the unit interval $[0,1]$ by a polygonal $p_1$ which has $N$ segments, all of them with length equal to $1/n$ ($N > n$) —that is, a replacement process of $n$ equal segments by $N$ equal segments. Successively, the polygonals $p_2, p_3, ...$, etc., are obtained by making the same $(n, N)$ substitution on each segment of the preceding polygonal.

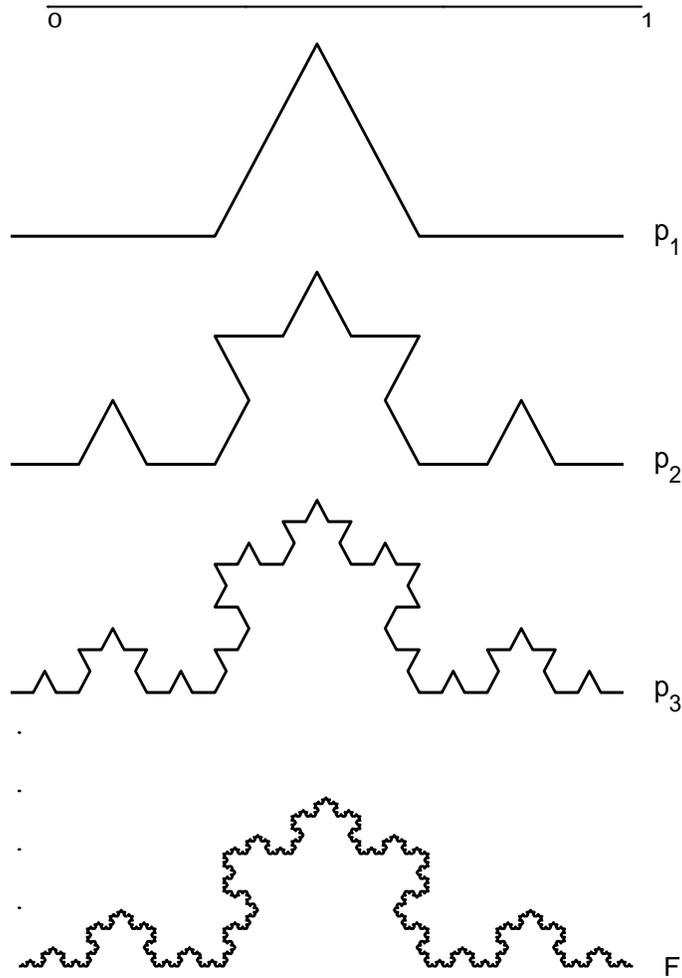

**Figure 3**

Notice that in this process we have $N$ contractive transformations of ratios $a_i = 1/n$, $i = 1, ..., N$. Repeating this replacement process *ad infinitum*, we obtain a bounded continuous curve $F$ of infinite length and infinite "creasing" that belongs to $\mathcal{F}_H$ and whose Hausdorff dimension is:

$$\dim_H(F) = \frac{\log N}{\log n} \quad .$$

Now, if we start again with the interval [0,1], but in the first step we construct a polygonal $p'_1$ with $N$ unit segments, and whose shape is the same as $p_1$, then the diameter of $p'_1$ will be $n$ times larger than the diameter of $p_1$ —$p'_1$ will be $p_1$ expanded by a ratio of $n$ to 1, $n$ is the inverse of the unique contractive factor involved in the fractal construction. In the second step, we construct a polygonal $p'_2$ identical to $p_2$ but with diameter $n^2$ times larger than that of $p_2$, and so on. In this way, we obtain a continuous unbounded curve $\Gamma$, locally rectifiable, that belongs to $\mathcal{F}_{MF}$; we will call $\Gamma$ strictly self-similar. We associate $F$ with $\Gamma$, and we are going to compare the dimensions $\dim_H(F)$ and $\dim_{MF}(\Gamma)$.

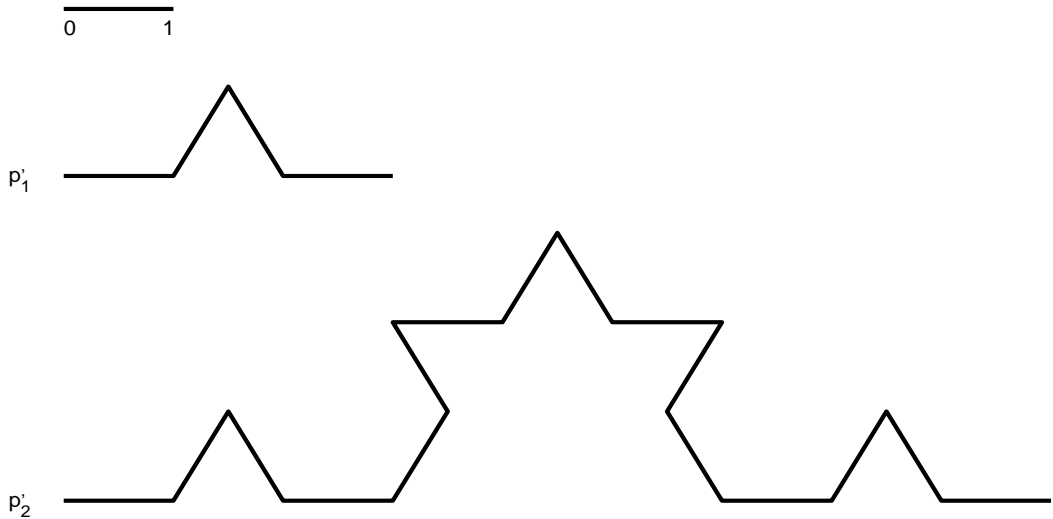

**Figure 4**

For any step $k$, the segments of the polygonal $p'_k$ are unity. Let $\ell_k$ be the length of $p'_k$, then, $A(\Gamma_k(\varepsilon)) \approx \varepsilon \times \ell_k = \varepsilon \times N^k$. If $C_k$ is the length of the boundary of the convex hull of $p'_k$, then $C_k \approx \text{const.} \times \text{diam}(p'_k) \approx \text{const.} \times n^k$. Thus, the dimension of Mendès France of the curve $\Gamma$ is:

$$\dim_{MF}(\Gamma) \;=\; \lim_{k\to\infty} \frac{\log A(\Gamma_k(\varepsilon))}{\log C_k} \;=\; \lim_{k\to\infty} \frac{\log(\varepsilon \times N^k)}{\log(\text{const.} \times n^k)} \;=\; \frac{\log N}{\log n} \;.$$

As we can note, in the strict self-similar case, the dimension of Mendès France of $\Gamma$ is equal to the Hausdorff dimension of the corresponding fractal $F$.

One question that arises in a natural way: does equality hold for processes other than the $(n, N)$ ones? If this equality does not hold, which, then, would be the relationship between these two dimensions: if, for instance, we take away strict self-similarity and allow $N$ contractive transformations whose ratios $a_i < 1$, $i = 1...N$ have **different** values?

4. SELF-SIMILAR CURVES

In the case of strict self-similarity, the expansive ratio used in the construction of polygonals $p'_k$ is the inverse of the **unique** contractive factor of the $N$ transformations generating the fractal. In this section we will study a more general type of curve, since we will allow the $N$ contractive ratios $a_i$ to be different. For example:

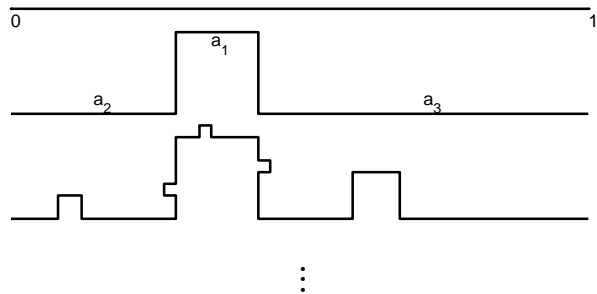

**Figure 5**: $N = 3$, $a_1 < a_2 < a_3 < 1$.

So, now, we have the possibility to construct a curve choosing an expansion ratio among the reciprocals $1/a_1, 1/a_2, ..., 1/a_N$, and to obtain $N$ curves $\Gamma^{a_1}, \Gamma^{a_2}, ..., \Gamma^{a_N}$, all of them the limit curve of a sequence of polygonals, all of them in $\mathcal{F}_{MF}$.

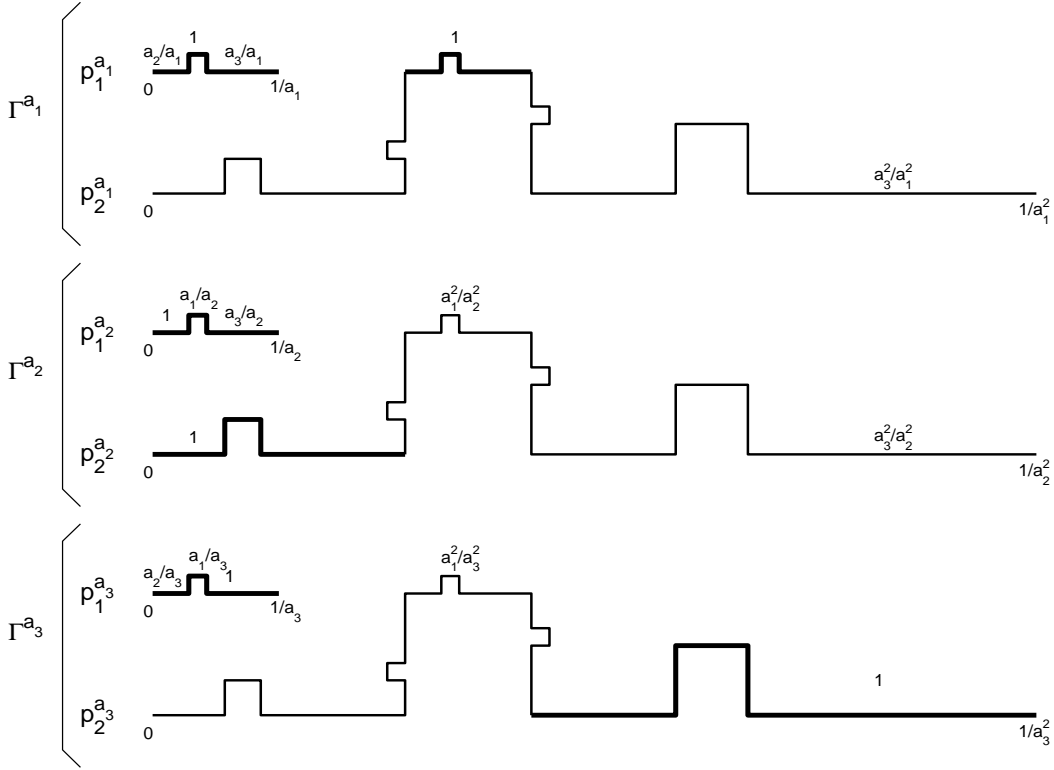

**Figure 6**: $N = 3$, $1 > a_1 > 1/a_2 > 1/a_3$.

We will analyze the geometric differences among $\Gamma^{a_1}, \Gamma^{a_2}, ..., \Gamma^{a_N}$; we will compare the different $d_{MF}$, and we will compare the latter with the **unique** Hausdorff dimension of the corresponding fractal $F$.

We can see (Fig. 6), in the first case, that the shortest segment is unity, and the others increase their length in each iteration. In the second case, the largest segment is unity, while the other segments have lengths going to zero in successive iterations. In the third situation, each iteration produces polygonal curves that have both larger and shorter segments than those in the preceding polygonal. This fact causes the three types of curve to be very different from one another. The first type is said to be a *resolvable curve*, while the others are *non-resolvable curves*.

The idea of the meaning of a resolvable curve is the following: We know that, for any curve $\Gamma \epsilon \mathcal{F}_{MF}$, if we take a closed ball in $I\!R^2$ with center on $\Gamma$, and we run this ball along the curve, the ball always contains a finite arc of $\Gamma$, but if this arc increases its length as the ball runs, i.e., if the arc into the ball is more "creased" and larger, then we say that $\Gamma$ is non-resolvable. Otherwise $\Gamma$ is resolvable (the formal definition is in [1]).

To sum up, let us consider an iterative-replacing-system that generates a fractal with ratios $a_1 \leq a_2 \leq ... \leq a_N < 1$; then each of its reciprocals $1/a_1 \geq 1/a_2 \geq ... \geq 1/a_N$ produces a different curve. If we take $1/a_1$ —the largest factor— the lengths of the segments will be always larger or equal to 1, and increasing in each iteration. If we take $1/a_N$ —the smallest factor— the length of segments in the successive polygonals decrease. If we take an intermediate factor $1/a_i$, $i \neq 1, N$, the polygonals will have both increasing and decreasing segments. Among those curves the only one that is resolvable is $\Gamma^{a_1}$, the first curve.

For this case, it is easy to calculate the dimension of Mendès France, because we can have a very good approximation of $A(\Gamma_k(\varepsilon))$ as $\varepsilon \times \ell^k$. Then, for all resolvable curves $\Gamma$, we have:

$$\dim_{MF}(\Gamma) = \lim_{k \to \infty} \frac{\log \ell_k}{\log C_k} \ .$$

If we bear in mind that in each step $k$, the length of polygonal $p'_k$ is:

$$\ell_k = \left(\frac{1}{a_1}\right)^k \left(\sum_{i=1}^{N} a_i\right)^k = \left(\frac{\sum_{i=1}^{N} a_i}{a_1}\right)^k$$

and

$$C_k \approx \text{const.} \times \left(\frac{1}{a_1}\right)^k$$

we then have

$$\dim_{MF}(\Gamma^{a_1}) = \frac{\log\left(\frac{\sum_{i=1}^{N} a_i}{a_1}\right)}{\log\left(\frac{1}{a_1}\right)} \ . \tag{1}$$

If we want to calculate the dimension of the curves in the other non-resolvable cases, it is very difficult to estimate the $A(\Gamma_k(\varepsilon))$ in a general manner. Nevertheless, if we call $\Gamma^{a_1}, \Gamma^{a_2}, ..., \Gamma^{a_N}$ the different curves, we are able to affirm that $\Gamma^{a_1}$ is the one with the minimal dimension of Mendès France, $\Gamma^{a_N}$ is that with the maximal dimension, and the rest of them have intermediate dimensions. The larger the dimension, the smaller the expansion factor.

This is the first part of Theorem 1, which will be proven later.

**Theorem 1- Part a).** If the contractive factors are $a_1 \leq a_2 \leq ... \leq a_N < 1$ and the expansive factors $1/a_1 \geq 1/a_2 \geq ... \geq 1/a_N$, then

$$\dim_{MF}(\Gamma^{a_1}) \leq \dim_{MF}(\Gamma^{a_2}) \leq ... \leq \dim_{MF}(\Gamma^{a_N}) \ .$$

How does the Hausdorff dimension of the fractal curve $F$ —built with the $a_i$ contractive factors— compare with all these different $d_{MF}$?

In order to tackle this subject, we will study the geometric differences of the expansive curves $\Gamma^{a_i}$, $i = 1, ..., N$. Consider one of these curves $\Gamma^{a_i}$, $i \neq N$. Then, in accordance with what we said above, it has segments as large as we want. Now, to fix ideas, let us consider $\varepsilon = 1$ and the 1-Minkowski sausage. Let us suppose that on each segment of this curve we make an *ad-infinitum* iteration of the corresponding process that generates the fractal $F$. This new fractal curve just doesn't belong to $\mathcal{F}_{MF}$, so there isn't a dimension of Mendès France associated with it. The curve doesn't belong to $\mathcal{F}_H$ either; however it has the same Hausdorff dimension $d$ as the fractal $F$, i.e., $d = \dim_H(F)$ —because there is a fractal like $F$ on each segment of $\Gamma^{a_i}$. The important thing to note here is that this curve is not covered by the 1-Minkowski sausage, since, if we take some segment in $\Gamma^{a_i}$ with length $\ell$ very very large, then the fractal built on it is not covered any more by the rectangle with area equal to $1 \times \ell$. And this fact is true **for all values of** $\varepsilon > 0$.

Let us take, now, the curve $\Gamma^{a_N}$ expanded by $1/a_N$. All the segments have length smaller than or equal to unity. Let us consider again the 1-Minkowski sausage of $\Gamma^{a_N}$ and let us make the same *ad-infinitum* iteration of the corresponding fractal $F$ on each segment of $\Gamma^{a_N}$. Thinking in the same manner as before, this new curve has a Hausdorff dimension equal to $\dim_H(F)$, but now this curve is completely covered by the 1-Minkowski sausage of $\Gamma^{a_N}$. If we took some $\varepsilon < 1$ in place of $\varepsilon = 1$, the curve will also be covered by the $\varepsilon$-sausage of $\Gamma^{a_N}$, except for the polygonals of the first iterations.

This means that $\Gamma^{a_N}$ is **the only curve** of all expanded curves that "shares" the $\Gamma(\varepsilon)$ with $F$ for every $\varepsilon$ "as if $\varepsilon$ were not be able to distinguish between a segment and a fractal built on it".

Since $\Gamma^{a_N}$ is the only curve, among all of $\Gamma^{a_i}$, that "gets creased" while its convex hull increases, this curve is the only one whose shape becomes the fractal form of $F$, as we move away from the plane.

This fact suggests that both dimensions, $\dim_H(F)$ and $\dim_{MF}(\Gamma^{a_N})$ are equal. This is the second part of Theorem 1.

**Theorem 1- Part b)** Let $F \epsilon \mathcal{F}_H$ be a fractal constructed by the contractive factors $a_1 \leq a_2 \leq ... \leq a_N < 1$, and let $\Gamma^{a_N}$ be the limit curve constructed by the expansive factor $1/a_N$. Then, we have:
$$\dim_H(F) = \dim_{MF}(\Gamma^{a_N}) \ .$$

Now, we are able to give the proof of Theorem 1.

### 4.1 PROOF OF THEOREM 1

**Proof of Part a).** Let $i$ be a fixed value, $1 \leq i \leq N$; let us suppose that $a_i = a_{i+1}^m$, where $m$ is an integer. Let us consider the corresponding polygonal $p_k$ of $\Gamma^{a_i}$ for a step $k$. Its length is:
$$\ell_k = \left( \frac{\sum_{j=1}^N a_j}{a_i} \right)^k = \left( \frac{\ell_0}{a_i} \right)^k \ ,$$

where $\ell_0 = \sum_{j=1}^n a_j$. The diameter in the $k$-step is
$$\operatorname{diam}(p_k) \approx \left( \frac{1}{a_i} \right)^k .$$

Let us now consider the step $m.k$. The corresponding polygonal of $\Gamma^{a_{i+1}}$, $p'_{mk}$, has length
$$\ell'_{mk} = \left( \frac{\ell_0}{a_{i+1}} \right)^{mk} \ ,$$

and diameter
$$\operatorname{diam}(p'_{mk}) \approx \left( \frac{1}{a_{i+1}} \right)^{mk} .$$

Then:
$$\ell'_{mk} = \frac{\ell_0^{mk}}{a_{i+1}^{mk}} = \frac{\ell_0^k}{a_i^k} \ell_0^{(m-1)k} = \ell_k \ell_0^{(m-1)k} \ ,$$

and then we have
$$\ell'_{mk} > \ell_k \ .$$

Therefore, for a fixed value of $\varepsilon$, we have
$$A(\Gamma_k^{a_i}(\varepsilon)) < A(\Gamma_{mk}^{a_{i+1}}(\varepsilon)) \ ,$$

and also
$$\operatorname{diam}(p'_{mk}) = \operatorname{diam}(p_k) \ ;$$

so we have
$$\frac{1}{a_{i+1}^{mk}} = \left( \frac{1}{a_{i+1}^m} \right)^k = \left( \frac{1}{a_i} \right)^k .$$

Thus
$$\frac{\log A(\Gamma_k^{a_i}(\varepsilon))}{\log \left( \frac{1}{a_i} \right)^k} < \frac{\log A(\Gamma_{mk}^{a_{i+1}}(\varepsilon))}{\log \left( \frac{1}{a_{i+1}} \right)^{mk}}$$

and taking limits when $k$ tends to infinity, we have:

$$\dim_{MF}(\Gamma^{a_i}) \leq \dim_{MF}(\Gamma^{a_{i+1}}) .$$

In the case $a_i = a_{i+1}^m$ with non integer $m$, the calculation is the same, considering the step $[m.k]$.

**Proof of Part b).** It is known that the Hausdorff dimension of a fractal $F$ with these characteristics can be expressed by:

$$\dim_H(F) = \lim_{\varepsilon \to 0}\left(2 - \frac{\log A(F(\varepsilon))}{\log(\varepsilon)}\right) , \qquad (2)$$

$A(F(\varepsilon))$ being the area of the $\varepsilon$-Minkowski sausage of $F$ ([2]). Let $(p_k)_{k \in I\!N}$ be the polygonals of successive $k$-steps in the fractal iteration of $F$, and let $(p'_k)_{k \in I\!N}$ be the corresponding polygonals of successive steps in the iteration of the limit curve $\Gamma$. For a certain value of $k$, we have that $p_k$ and $p'_k$ are "alike", that is, they have identical shape but different size. Then, if $\ell_k = $ length of $p_k$ and $\ell'_k = $ length of $p'_k$, we have

$$\ell'_k = \left(\frac{1}{a_N}\right)^k \ell_k .$$

Besides, if we consider the area of the $\varepsilon$-sausage of $p'_k$, we have that the corresponding "alike" sausage of $p_k$ satisfies:

$$A(p'_k(\varepsilon)) = \left(\frac{1}{a_N}\right)^{2k} A(p_k(\varepsilon \times a_N^k)) . \qquad (3)$$

On the other hand, we have that for every polygonal $p_k$, the largest segment has length equal to $a_N^k$; so, for every $k$ we can write:

$$A(p_k(a_N^k)) \approx A(F(a_N^k)) . \qquad (4)$$

Notice that $A(p_k(\lambda^k)) \approx A(F(\lambda^k))$ is not valid for all $k$ when $\lambda < a_N$; for, as $k$ grows, so does the difference between $\lambda^k$ and $a_N^k$, breaking down the comparability stated in (4). Therefore, taking $\varepsilon = 1$, from (3) and (4) we have:

$$A(p'_k(1)) \approx \left(\frac{1}{a_N}\right)^{2k} A(F(a_N^k)) .$$

Then:

$$\dim_{MF}(\Gamma) = \lim_{k \to \infty} \frac{\log A(p'_k(1))}{\log\left(\frac{1}{a_N}\right)^k} = \lim_{k \to \infty} \frac{\log\left(\left(\frac{1}{a_N}\right)^{2k} A(F(a_N^k))\right)}{\log\left(\frac{1}{a_N}\right)^k} = \lim_{k \to \infty}\left(2 - \frac{\log A(F(a_N^k))}{\log(a_N^k)}\right)$$

and then, by (2):

$$\dim_{MF}(\Gamma) = \dim_H(F) .$$

The following theorem states that the inequalities from the first part of Theorem 1 are strict inequalities.

4.2. THEOREM 2

**Theorem 2.** Let $\ell, i \in I\!N$ be such that $\ell < i$ and $a_\ell < a_i$; let $\Gamma^{a_\ell}$ and $\Gamma^{a_i}$ be the expanded limit curves corresponding to the factors $1/a_\ell$ and $1/a_i$ respectively. Then the following inequality holds:

$$\dim_{MF}(\Gamma^{a_\ell}) < \dim_{MF}(\Gamma^{a_i}) .$$

**Proof.** We will make here a simplification, considering the case $\ell = 1$ and $1 < i \leq N$, since the proof for the general case —such as it is at present— would exceed the limits of this work.

First we need two results:

1) Let $\alpha, k \in \mathbb{N}, k \geq 2$, then

$$\sum_{i=0}^{k/2} \binom{k}{i} \left(\frac{1}{\alpha}\right)^i \geq \frac{1}{2}\left(1 + \frac{1}{\alpha}\right)^k \quad ,$$

and if $k$ is odd, we replace $k/2$ by $(k-1)/2$ —so it is sufficient to consider $k$ even.

In fact,

$$\binom{k}{0}\left(\frac{1}{\alpha}\right)^0 \geq \binom{k}{k}\left(\frac{1}{\alpha}\right)^k$$

$$\binom{k}{1}\left(\frac{1}{\alpha}\right)^1 \geq \binom{k}{k-1}\left(\frac{1}{\alpha}\right)^{k-1}$$

$$\vdots$$

$$\binom{k}{\frac{k}{2}-1}\left(\frac{1}{\alpha}\right)^{\frac{k}{2}-1} \geq \binom{k}{\frac{k}{2}+1}\left(\frac{1}{\alpha}\right)^{\frac{k}{2}+1}.$$

Then,

$$\sum_{i=0}^{\frac{k}{2}-1} \binom{k}{i}\left(\frac{1}{\alpha}\right)^i \geq \sum_{i=\frac{k}{2}+1}^{k} \binom{k}{i}\left(\frac{1}{\alpha}\right)^i \quad ,$$

and therefore

$$\sum_{i=0}^{\frac{k}{2}} \binom{k}{i}\left(\frac{1}{\alpha}\right)^i > \sum_{i=\frac{k}{2}+1}^{k} \binom{k}{i}\left(\frac{1}{\alpha}\right)^i.$$

To both sides of the last inequality we add the first side, and we obtain:

$$2\sum_{i=0}^{\frac{k}{2}} \binom{k}{i}\left(\frac{1}{\alpha}\right)^i > \sum_{i=0}^{k} \binom{k}{i}\left(\frac{1}{\alpha}\right)^i \quad ,$$

that is to say

$$\sum_{i=0}^{\frac{k}{2}} \binom{k}{i}\left(\frac{1}{\alpha}\right)^i > \frac{1}{2}\left(1 + \frac{1}{\alpha}\right)^k.$$

2) If we take the $k$th-step in the construction of the limit curve $\Gamma^{a_i}$ — $i$ fixed, expansion factor $1/a_i$— we see that the segments of the polygonal $p'_k$ have lengths $\frac{A}{B}$; the denominator $B$ is always $a_i^k$, the numerator $A$ is always a product of powers of $a_1, a_2, ..., a_N$ in such a way that the sum of their exponents is equal to $k$:

$$A = a_1^{j_1} a_2^{j_2} ... a_i^{j_i} ... a_N^{j_N} \quad , \qquad \sum_{\alpha=1}^{N} j_\alpha = k \quad .$$

For one such configuration $(j_1, j_2, ..., j_N)$ we pose the question: how many segments of **this** length are there in step $k$? Answer: the number of such segments is the numerical coefficient of the term whose "literal" part is $a_1^{j_1} a_2^{j_2} ... a_i^{j_i} ... a_N^{j_N}$, in the development of $(a_1 + a_2 + ... + a_N)^k$. That is to say, there are:

$$\binom{k}{j_1}\binom{k-j_1}{j_2}\binom{k-j_1-j_2}{j_3}...\binom{k-j_1-j_2-...-j_{N-2}}{j_{N-1}}\binom{k-j_1-j_2-...-j_{N-1}}{j_N}.$$

Let us now consider those segments for which the numerator $A$ has the following configuration: half —or less— of the factors that appear in $A$ are equal to $a_1$, and the rest of them —to complete a total of $k$ factors— are permutations of $a_2, a_3, ..., a_N$. That is to say, $j_1 \in [0, k/2]$.

Now we will prove that the quantity of those segments is larger than or equal to the half of $N^k$, which is the total number of segments appearing in step $k$.

Indeed, let us fix exponent $j_\alpha$, $\alpha \neq N-1$, $j_{N-1}$ will run from 0 to $k - j_1 - j_2 - ... - j_{N-2}$; and let us count the corresponding number of such segments:

$$\sum_{j_{N-1}=0}^{k-j_1-j_2-...-j_{N-2}} \binom{k}{j_1}\binom{k-j_1}{j_2}...\binom{k-j_1-j_2-...-j_{N-2}}{j_{N-1}} =$$

$$= \binom{k}{j_1}\binom{k-j_1}{j_2}...\binom{k-j_1-j_2-...-j_{N-3}}{j_{N-2}} \sum_{j_{N-1}=0}^{k-j_1-j_2-...-j_{N-2}} \binom{k-j_1-j_2-...-j_{N-2}}{j_{N-1}} =$$

$$= \binom{k}{j_1}\binom{k-j_1}{j_2}...\binom{k-j_1-j_2-...-j_{N-3}}{j_{N-2}} 2^{k-j_1-j_2-...-j_{N-2}}.$$

Now, let us fix all $j_\alpha$, except $j_{N-2}$ which will be allowed to run from 0 to $k-j_1-j_2-...-j_{N-3}$, and let us obtain the total of the corresponding segments:

$$\sum_{j_{N-2}=0}^{k-j_1-j_2-...-j_{N-3}} \binom{k}{j_1}\binom{k-j_1}{j_2}...\binom{k-j_1-j_2-...-j_{N-3}}{j_{N-2}} 2^{k-j_1-j_2-...-j_{N-3}} 2^{-j_{N-2}} =$$

$$= 2^{k-j_1-j_2-...-j_{N-3}} \binom{k}{j_1}...\binom{k-j_1-j_2-...-j_{N-4}}{j_{N-3}} \sum_{j_{N-2}=0}^{k-j_1-j_2-...-j_{N-3}} \binom{k-j_1-j_2-...-j_{N-3}}{j_{N-2}} 2^{-j_{N-2}} =$$

$$= 2^{k-j_1-j_2-...-j_{N-3}} \binom{k}{j_1}...\binom{k-j_1-j_2-...-j_{N-4}}{j_{N-3}} \left(1+\frac{1}{2}\right)^{k-j_1-j_2-...-j_{N-3}}.$$

Next, let us fix all $j_\alpha$, except $j_{N-3}$, that we will run from 0 to $k - j_1 - j_2 - ... - j_{N-4}$, doing the same as before, and we again obtain the total number of corresponding segments:

$$\binom{k}{j_1}\binom{k-j_1}{j_2}...\binom{k-j_1-j_2-...-j_{N-5}}{j_{N-4}} 4^{k-j_1-...-j_{N-4}}.$$

Proceeding in this way with the rest of the $j_\alpha$, let us fix $j_1$ and let us run $j_2$ from 0 to $k - j_1$. Counting as above, the corresponding number of such segments is:

$$\binom{k}{j_1}(N-1)^{k-j_1}, \qquad (j_1 = j_{N-(N-1)}).$$

Finally, let us run $j_1$ from 0 to $k/2$. We obtain the total number of segments:

$$\sum_{j_1=0}^{k/2} \binom{k}{j_1}(N-1)^{k-j_1} = (N-1)^k \sum_{j_1=0}^{k/2} \binom{k}{j_1}(N-1)^{-j_1} \geq$$

$$\geq (N-1)^k \frac{1}{2}\left(1+\frac{1}{N-1}\right)^k = \frac{1}{2}N^k,$$

using our result 1).

We can prove the theorem now. We want to prove that if $i > 1$, then $\dim_{MF}(\Gamma^{a_i}) > \dim_{MF}(\Gamma^{a_1})$. We will suppose, without loss of generality, that $a_1 < a_2 \leq a_i$, $i$ fixed between 2 and $N$.

Let $p_k^i$ be the polygonal of the corresponding $k$th step of the limit curve $\Gamma^{a_i}$. The diameter of this polygonal is equal to $\left(\frac{1}{a_i}\right)^k$.

As we said before, in this polygonal there are $N^k$ segments whose lengths can be written thus:

$$\frac{a_1^{j_1} a_2^{j_2} ... a_i^{j_i} ... a_N^{j_N}}{a_i^k} \quad , \quad \sum_{\alpha=1}^{N} j_\alpha = k \ .$$

Let $m \in \mathbb{N}, m < k$, and let $p_m^i$ be the corresponding polygonal. This polygonal has $N^m$ segments, and a diameter equal to $\left(\frac{1}{a_i}\right)^m$.

Let $\overline{p}_m^i$ be a polygonal with the same shape as $p_m^i$ but expanded by a ratio equal to $\left(\frac{1}{a_i}\right)^{k-m}$.

The diameter of $\overline{p}_m^i$ is, now, $\left(\frac{1}{a_i}\right)^k$, and the lengths of its segments are of the form:

$$\left(\frac{1}{a_i}\right)^{k-m} \frac{a_1^{j_1} a_2^{j_2} ... a_N^{m-(j_1+...+j_{N-1})}}{a_i^m} \quad ,$$

where $j_\alpha$ are not necessarily the same as before.

Now let us choose a number $m$ —later we will exhibit the explicit value of this $m$— such that any segment with length equal to

$$\left(\frac{1}{a_i}\right)^{k-m} \frac{a_1^{m/2} a_2^{m/2} a_3^0 ... a_N^0}{a_i^m} = \frac{a_1^{m/2} a_2^{m/2}}{a_i^k}$$

—that is to say $j_1 = j_2 = m/2$ and $j_3 = ... = j_N = 0$— is comparable to unity.

This entails that segments, whose lengths have a "configuration" such that $j_1 \leq m/2$, the remaining $j_\alpha$ arbitrary —the sum always being $m$— become larger than or equal to unity. Now then, the number of these segments is, by our result 2), larger than or equal to the half of the total number of all segments, which is $N^m$.

In other words, if $\ell_m$ is the length of $p_m^i$, and $\overline{\ell}_m$ is the length of $\overline{p}_m^i$, we have:

$$\overline{\ell}_m = \left(\frac{1}{a_i}\right)^{k-m} \quad , \quad \ell_m = \left(\frac{1}{a_i}\right)^{k-m} \frac{\left(\sum_{j=1}^N a_j\right)^m}{a_i^m} = \frac{\left(\sum_{j=1}^N a_j\right)^m}{a_i^k} \quad ,$$

and also:

$$\overline{\ell}_m < \ell_k \ .$$

Therefore, taking $\varepsilon = 1$, it is true that:

$$A(\overline{p}_m^i(1)) \leq A(p_k^i(1)) \ .$$

But, since more than half of the total number of segments are segments larger than or equal to unity, it follows that:

$$1 \times \frac{1}{2}\overline{\ell}_m \leq A(\overline{p}_m^i(1))$$

and then

$$\frac{1}{2}\frac{\left(\sum_{j=1}^N a_j\right)^m}{a_i^k} \leq A(p_k^i(1)) \ .$$

Therefore,
$$\frac{\log\left(\frac{1}{2}\frac{\left(\sum_{j=1}^{N} a_j\right)^m}{a_i^k}\right)}{\log\left(\frac{1}{a_i^k}\right)} \leq \frac{\log A(p_k^i(1))}{\log\left(\frac{1}{a_i^k}\right)} .$$

Next, we will calculate the explicit value of $m$ —which we have chosen in order to satisfy the last inequality. The value $m$ was chosen requiring that:
$$\frac{a_1^{m/2} a_2^{m/2}}{a_i^k} \approx 1 ;$$

then,
$$\frac{a_1^{m/2} a_2^{m/2}}{a_i^k} = \frac{a_1^{m/2}\left(a_1^{\log_{a_1}(a_2)}\right)^{m/2}}{a_i^k} =$$
$$= \frac{a_1^{m/2(\log_{a_1}(a_2)+1)}}{a_i^k} = \frac{\left(a_1^{\frac{m(1+\log_{a_1}(a_2))}{2k}}\right)^k}{a_i^k} \approx 1 .$$

That is to say:
$$\frac{a_1^{m(1+\log_{a_1}(a_2))}}{2k} \approx a_i ;$$

therefore
$$\frac{m(1+\log_{a_1}(a_2))}{2k} \approx \frac{1}{\log_{a_i}(a_1)} ,$$

that is to say
$$\frac{m}{k} \approx \frac{2}{\log_{a_i}(a_1)(1+\log_{a_1}(a_2))} = \frac{2}{\log_{a_i}(a_1) + \log_{a_i}(a_1)\log_{a_1}(a_2)} ,$$

so that
$$\frac{m}{k} \approx \frac{2}{\log_{a_i}(a_1 a_2)} .$$

Going back to inequality (5), we have that:
$$\frac{\log\left(\frac{1}{2}\right) + m \log\left(\sum_{j=1}^{N} a_j\right) + k \log\left(\frac{1}{a_i}\right)}{k \log\left(\frac{1}{a_i}\right)} \leq \frac{\log A(p_k^i(1))}{l \log\left(\frac{1}{a_i}\right)} .$$

Replacing $\frac{m}{k}$ by the expression (6), we obtain:
$$\frac{\log\left(\frac{1}{2}\right)}{k \log\left(\frac{1}{a_i}\right)} + \frac{2}{\log_{a_i}(a_1 a_2)} \frac{\log\left(\sum_{j=1}^{N} a_j\right)}{\log\left(\frac{1}{a_i}\right)} + 1 \leq \frac{\log A(p_k^i(1))}{l \log\left(\frac{1}{a_i}\right)} ,$$

and taking limits when $k$ tends to infinity, we have that:
$$\frac{\log\left(\sum_{j=1}^{N} a_j\right)}{\log\left(\frac{1}{\sqrt{a_1 a_2}}\right)} + 1 \leq \dim_{MF}(\Gamma^{a_i}) .$$

Finally, since $a_1 < a_2$, we have $\frac{1}{a_1} > \frac{1}{\sqrt{a_1 a_2}}$, so that:
$$1 + \frac{\log\left(\sum_{j=1}^{N} a_j\right)}{\log\left(\frac{1}{a_1}\right)} < 1 + \frac{\log\left(\sum_{j=1}^{N} a_j\right)}{\log\left(\frac{1}{\sqrt{a_1 a_2}}\right)} .$$

Therefore, taking inequality (1) into account, we have:

$$\dim_{MF}(\Gamma^{a_1}) < \dim_{MF}(\Gamma^{a_i}) .$$

q.e.d.

5. REMARK

The preceding theorem's proof is based on arguments and ideas that are strongly geometrical; however, if we consider the extremal cases $\ell = 1$ and $i = N$ —minimal and maximal dimensions respectively— a completely different —and much shorter!— proof of the corresponding inequality:

$$\dim_{MF}(\Gamma^{a_1}) < \dim_{MF}(\Gamma^{a_N}) , \qquad (5)$$

can be given. Indeed, for the first curve $\Gamma^{a_1}$ we have a formula that allows us to calculate its dimension, and for $\Gamma^{a_N}$ —which is the only curve whose dimension is equal to the dimension of the associated fractal, by virtue of Theorem 1, Part b)— we have an implicit equation that is satisfied by the Hausdorff dimension of the fractal.

**Proof of inequality (5).** The Hausdorff dimension is the value $d_H$ that satisfies:

$$\sum_{i=1}^{N} a_i^{d_H} = 1 ;$$

on the other hand, if $d = \dim_{MF}(\Gamma^{a_1})$, then $d$ satisfies:

$$\sum_{i=1}^{N} \frac{a_i}{a_1} a_1^d = 1 ,$$

since we have

$$\sum_{i=1}^{N} a_i a_1^{d-1} = 1 ,$$

if and only if

$$a_1^d \sum_{i=1}^{N} \frac{a_i}{a_1} = 1 .$$

Taking logarithms in this last equation, we obtain:

$$d \, \log(a_1) + \log\left(\frac{\sum_{i=1}^{N} a_i}{a_1}\right) = 0 ,$$

and then:

$$d = \frac{\log\left(\frac{\sum_{i=1}^{N} a_i}{a_1}\right)}{\log\left(\frac{1}{a_1}\right)} = \dim_{MF}(\Gamma^{a_1}) .$$

Let us consider now the following functions:

$$f(x) = \left(\sum_{i=1}^{N} a_i a_1^{x-1}\right) - 1 \quad \text{and} \quad g(x) = \left(\sum_{i=1}^{N} a_i^x\right) - 1 .$$

According to what we just wrote, we have:

$$f(d) = 0 \quad \text{and} \quad g(d_H) = 0 .$$

Besides, $f$ and $g$ are decreasing functions ($a_i < 1, i = 1, ..., N$). Let us compare any term in the expression of $f$ with the corresponding term in $g$ —except the first term which is equal in both functions:
for $x > 1$, we have
$$a_1^{x-1} \;<\; a_i^{x-1} \quad (i \neq 1) \;,$$
so
$$\frac{a_i}{a_1} a_1^x \;\leq\; a_i^x \;;$$
therefore
$$f(x) \;<\; g(x) \quad \text{for } x > 1 \;;$$
in particular
$$f(d_H) \;<\; g(d_H) \;,$$
and therefore
$$d \;<\; d_H \;.$$

NOTE. From this remark we conclude that the study of Hausdorff and Mendès France dimensions of curves $\Gamma^{a_i}$, $1 \leq i \leq N$, associated with the same fractal curve $F$, has a very different nature for the case $i = 1$, $i = N$, and for the case $i \neq 1, N$. Because of this we gave the proof of Theorem 2 restricting the general situation $1 \leq \ell < i \leq N$ to $\ell = 1$, since, as stated above, the proof of the general case, such as it is at present, would exceed the limits of this paper.

ACKNOWLEDGMENT. The authors wish to thank E. Cesaratto for the valuable help in the preparation of all the figures in this article.